\documentclass[11pt]{article}
\usepackage[T1]{fontenc}
\usepackage{lmodern}
\usepackage[a4paper,left=1.02in,right=1.02in,top=.92in,bottom=.92in]{geometry}
\usepackage{amsmath,amssymb,booktabs}
\usepackage[hidelinks]{hyperref}
\hypersetup{pdftitle={Two-weight enumeration and uniform ratio limits for closed meanders},pdfauthor={Iaroslav Argunov}}
\newtheorem{theorem}{Theorem}
\newtheorem{lemma}[theorem]{Lemma}
\newtheorem{corollary}[theorem]{Corollary}
\newenvironment{proof}{\par\noindent\textit{Proof.}\ }{\unskip\nobreak\hfill\(\square\)\par}
\newcommand{\M}{\mathcal M}
\newcommand{\Prb}{\mathbb P}
\newcommand{\E}{\mathbb E}
\newcommand{\bin}[2]{\binom{#1}{#2}}

\begin{document}
\begin{center}
{\LARGE\bfseries Two-weight enumeration and uniform ratio limits\\[.1em] for closed meanders\par}
\vspace{.55em}
{\large Iaroslav Argunov\par}
{\small Independent researcher, Oakland, USA\par}
{\small Corresponding author: \nolinkurl{argunov88@gmail.com}\par}
\vspace{.15em}
{November 23, 2026}
\end{center}
\begin{abstract}
In a closed meander, two consecutive arches on opposite sides of the line
can occur in either order. We count these patterns on triples beginning
at an even label, giving their two orientations separate weights \(u\)
and \(v\). Contracting runs of marked sites expresses the weighted
generating function as a rational substitution in the ordinary meander
series. We prove that the ratio of consecutive weighted counts converges
uniformly on each compact subset of the closed nonnegative quadrant.
The substitution determines the limit in terms of the ordinary growth
constant, including when both patterns are forbidden. At zero weights
the proof uses two coprime return lengths in a transfer that remembers
the preceding crossing move. For unit weights the ratio limit also
follows from Jensen's transfer and Kesten's ratio limit lemma.
\end{abstract}

\section{Two weights for adjacent arches}
A closed meander on \(2n\) points is an ordered pair of noncrossing
perfect matchings on \(0,1,\ldots,2n-1\) whose union is connected.
Draw the first matching above the line and the second below.
The order of the points and the choice of half-planes are fixed.
At \(n=1\) the union consists of two parallel edges. The first counts
are \(M_1=1,M_2=2,M_3=8,M_4=42\).

There are \(n-1\) test sites. Site \(i\), \(0\le i\le n-2\), is
\emph{positive} when the top matching joins \(2i\) to \(2i+1\) and
the bottom matching joins \(2i+1\) to \(2i+2\). A \emph{negative}
site uses the opposite colours. The line has two ends, so no site
wraps around. Write \(X_n\) for the number of positive sites in a
uniform connected meander.

\begin{center}
\setlength{\unitlength}{1pt}
\begin{picture}(206,49)
\put(1,26){\line(1,0){85}}\put(121,26){\line(1,0){85}}
\qbezier(14,26)(30,54)(46,26)\qbezier(46,26)(62,-2)(78,26)
\qbezier(134,26)(150,-2)(166,26)\qbezier(166,26)(182,54)(198,26)
\multiput(14,26)(32,0){3}{\circle*{3.4}}
\multiput(134,26)(32,0){3}{\circle*{3.4}}
\put(32,0){\makebox(30,8){\small positive}}\put(152,0){\makebox(30,8){\small negative}}
\end{picture}
\end{center}

Give each positive occurrence weight \(u\) and each negative occurrence
weight \(v\), and let \(Z_n(u,v)\) be the total weight of the order-\(n\)
meanders. The ordinary series is \(\M(z)=\sum_{n\ge1}M_nz^n\), with
growth constant \(\mu=\lim_{n\to\infty}M_n^{1/n}\) \cite{AP}. The two
weights enter through a single substitution (Theorem~\ref{thm:joint}):
\[
  \sum_{n\ge1}Z_n(u,v)z^n=\M(f_{u,v}(z)),\qquad
  f_{u,v}(z)=z\left(\frac{1}{1-(u-1)z}+\frac{1}{1-(v-1)z}-1\right).
\]
If \(\rho(u,v)\) is the first positive solution of
\(f_{u,v}(\rho)=1/\mu\) before any positive pole, then
\(Z_{n+1}(u,v)/Z_n(u,v)\to1/\rho(u,v)\) uniformly on every compact
subset of \([0,\infty)^2\)
(Theorem~\ref{thm:all-weighted-ratios} and
Corollary~\ref{cor:uniform-weighted-ratios}).
At \(u=v=0\) the two patterns are forbidden and \(Z_2(0,0)=0\),
but the counts are positive for \(n\ge3\). Their consecutive ratio tends to
\((\mu-1+\sqrt{\mu^2-6\mu+1})/2\).

Arch statistics in folding models were studied by Di Francesco, Golinelli
and Guitter \cite{DGG}. Here all diagrams are connected, and a prescribed
run at specified sites can be contracted to one vertex. This accounts
for the substitution. To pass from the radius of convergence to
consecutive coefficient ratios, we keep the last crossing move in the
transfer state. Two unmarked returns then survive at every choice of
weights. At \(u=v=1\) the ordinary ratio limit also follows by applying
Kesten's lemma \cite[Lemma~4]{K95} to Jensen's transfer \cite{Jensen}.

\section{Contraction and exchangeability}
Throughout, \([h]\) denotes a set of \(h\) consecutive sites, with
\([0]=\varnothing\). Write \(B_n\) for the number of order-\(n\)
connected meanders with no positive site.

\begin{lemma}\label{lem:contract}
Specify signs on a set \(S\subseteq[n-1]\) of \(k\) sites. The number
of meanders realising all prescriptions is \(M_{n-k}\), unless two
adjacent specified sites have opposite signs, in which case it is zero.
\end{lemma}
\begin{proof}
Adjacent sites share one vertex. Incompatible signs require it to have
two incident arches of the same colour. With compatible signs, a run of
\(s\) consecutive marked sites is an alternating path on \(2s+1\)
consecutive vertices. There are two free incidences, one at each end
and of opposite colours. Collapse the path to a single vertex with
these two incidences. This shortens the unique cycle and preserves the
ordering of the exterior arches in each half-plane. The runs occupy
disjoint intervals, and contraction removes \(2k\) vertices in total.

Given the smaller meander, the marked positions determine which vertices
to expand. The prescribed sign of each run determines where its upper
and lower exterior arches attach. Every new interior arch joins
consecutive labels, so none crosses an exterior arch. The operations
are inverses, including when the two exterior arches reach the same vertex.
\end{proof}

\begin{theorem}\label{thm:law}
For every \(S\subseteq[n-1]\) the number of meanders whose set of
positive sites is \emph{exactly} \(S\) equals \(B_{n-|S|}\). In
particular,
\begin{equation}\label{eq:bin}
 M_n=\sum_{k=0}^{n-1}\bin{n-1}{k}B_{n-k},\qquad
 \Prb(X_n=k)=\bin{n-1}{k}\frac{B_{n-k}}{M_n}.
\end{equation}
\end{theorem}
\begin{proof}
Excluding positive occurrences at all sites outside \(S\) costs, by
inclusion--exclusion and Lemma~\ref{lem:contract},
\[
  \sum_{h=0}^{n-1-|S|}(-1)^h\bin{n-1-|S|}{h}M_{n-|S|-h}.
\]
Taking \(S=\varnothing\) at order \(n-|S|\) identifies this as
\(B_{n-|S|}\). Summing over sets of size \(k\) gives the two formulas.
\end{proof}

There is also a direct inverse for this exact-set contraction. Along
a positive run the lower exterior incidence lies at its left end,
the upper at its right end. The tests for positivity at the two
neighbouring exterior sites therefore survive contraction unchanged.
Expanding a chosen even vertex of an avoiding diagram into a prescribed
positive run restores exactly its marked sites.

\section{Exponential rates and the ratio theorem}
The avoidance counts have their own growth constant. To see this, lay
two avoiding meanders of orders \(r,s\) side by side. Remove the upper
arch \((a,2r-1)\) at the right end of the first and the upper arch
\((2r,b)\) at the left end of the second. Insert
\((a,b)\) and \((2r-1,2r)\). The first new arch encloses the two
pieces formerly under the removed arches; the second spans consecutive
labels. The operation merges the two cycles, and their original ends
can be recovered at the fixed seam. The short new upper arch starts at
an odd label and the other has length at least three. Neither is the
upper arch of a positive site. A possible new triple at the seam would
need a lower arch crossing the seam or a new positive upper arch.
Hence \(B_{r+s}\ge B_rB_s\). The construction also works without the
avoidance condition.

Since \(B_1=1\), these inequalities and the Catalan-square upper bound
give finite limits \(b=\lim B_n^{1/n}\ge1\) and
\(\mu=\lim M_n^{1/n}\). For \(0<\epsilon<b\) there are positive
constants \(c,C\) such that
\(c(b-\epsilon)^j\le B_j\le C(b+\epsilon)^j\) for every \(j\ge1\).
Substituting these bounds in~\eqref{eq:bin} and taking \(n\)-th roots gives
\begin{equation}\label{eq:base}
                      \mu=1+b .
\end{equation}

\begin{theorem}\label{thm:ratio}
Put \(p=1/\mu\). For every sequence \(0\le k_n\le n-1\) with
\(k_n/(n-1)\to x\in[0,1]\),
\[
 \lim_{n\to\infty}\frac1n\log\Prb(X_n=k_n)
 =-D(x\Vert p),\qquad
 D(x\Vert p)=x\log\frac{x}{p}+(1-x)\log\frac{1-x}{1-p},
\]
with \(0\log0=0\). Furthermore, \(\displaystyle\lim_{n\to\infty}M_n/M_{n-1}=\mu\).
\end{theorem}
\begin{proof}
In~\eqref{eq:bin}, the entropy estimate for binomial coefficients gives
\(n^{-1}\log\bin{n-1}{k_n}\to-x\log x-(1-x)\log(1-x)\).
Also \(n^{-1}\log B_{n-k_n}\to(1-x)\log b\). At \(x=1\) the
last conclusion follows directly from \(1\le B_j\le 16^j\).
Subtract \(\log\mu\) and use~\eqref{eq:base}. This proves the pointwise
claim at both endpoints as well as in the interior.

The rate vanishes only at \(p\). On any closed subset of \([0,1]\)
that excludes \(p\), bound the probability by \(n\) times its largest
point mass. A convergent subsequence of maximising densities gives
exponential decay by the pointwise claim. Consequently
\(\E X_n/(n-1)\to p\). For each individual site,
Lemma~\ref{lem:contract} gives probability \(M_{n-1}/M_n\).
Taking the mean of the indicators proves the stated ratio.
\end{proof}

\subsection*{A second proof of the ratio limit}
The avoidance splice can be omitted for the ratio assertion alone.
In the following lemma take \(a_N=M_{N+1}\) and \(b_j=B_{j+1}\)
using \eqref{eq:bin}. The only growth input is the ordinary root limit.
\begin{lemma}\label{lem:binomial-smoothing}
Let \(b_j\geq 0\), \(b_0>0\), and
\(a_N=\sum_{j=0}^N\binom Nj b_j\).
If \(a_N^{1/N}\to\lambda\) with \(1<\lambda<\infty\), then
\(a_{N-1}/a_N\to 1/\lambda\).
\end{lemma}
\begin{proof}
Write \(A(z)=\sum a_Nz^N\) and \(B(t)=\sum b_jt^j\).
For \(0\leq z<1\), summation of nonnegative terms gives
\(A(z)=(1-z)^{-1}B(z/(1-z))\), allowing infinite sums.
Since \(A\) has radius \(1/\lambda\), the radius of \(B\) is
\(1/(\lambda-1)\). For every \(\epsilon>0\) we thus have
\(b_j\leq C_\epsilon(\lambda-1+\epsilon)^j\) for all \(j\).

Give \(J_N\) the law
\(\Pr(J_N=j)=\binom Nj b_j/a_N\) and put \(p=(\lambda-1)/\lambda\).
For \(t=j/N\), the binomial theorem gives \(\binom Nj\leq e^{NH(t)}\),
where \(H(t)=-t\log t-(1-t)\log(1-t)\), with \(0\log0=0\).
The continuous function
\[
 f(t)=H(t)+t\log(\lambda-1)-\log\lambda
\]
has unique maximum \(f(p)=0\) on \([0,1]\).
Fix \(\eta>0\) such that \(K=\{t\in[0,1]:|t-p|\geq\eta\}\) is
nonempty, and set \(d=-\max_{K}f>0\).
Choose \(\epsilon>0\) with
\(\log(1+\epsilon/(\lambda-1))<d/3\).
For sufficiently large \(N\),
\(N^{-1}\log a_N\geq\log\lambda-d/3\).
Consequently
\[
 \Pr(J_N/N\in K)\leq (N+1)C_\epsilon e^{-Nd/3}.
\]
It follows that \(\mathbb E[J_N/N]\to p\). The identity
\(\binom{N-1}{j}=(1-j/N)\binom Nj\) now gives
\(a_{N-1}/a_N=1-\mathbb E[J_N/N]\to1/\lambda\).
\end{proof}

\section{The two signs and longer intervals}
\begin{corollary}
For any fixed \(k\) distinct sites, their positive indicators tend
jointly to independent Bernoulli variables of parameter \(1/\mu\).
The sites may vary with \(n\). With both signs recorded, along any set
of sites with fixed adjacency pattern the limiting joint intersection
probability equals \(\mu^{-k}\) for compatible prescribed signs and
zero otherwise.
\end{corollary}
\begin{proof}
For every prescribed subset, apply Lemma~\ref{lem:contract} and
\(M_{n-k}/M_n\to\mu^{-k}\), obtained by telescoping Theorem~\ref{thm:ratio}.
Inclusion--exclusion gives each exact assignment, including vacant
sites. For positive indicators the intersection probabilities factor.
For both signs they factor between sets separated by at least one unused
site. These consistent finite laws define a stationary one-dependent
process on the site lattice. Adjacent opposite signs have probability
zero, so the coloured process is not independent.
\end{proof}

The signed limit also has an explicit sampling rule. Make two vertices
\((i,+),(i,-)\) at each site, with edges between opposite signs at
distances zero and one. Orient each edge independently and uniformly;
mark each vertex independently with probability \(8/\mu\), which is
less than one by \cite{AP}. Occupy a marked vertex if all three of
its edges point towards it. Adjacent vertices cannot both be occupied,
whereas any independent set of \(k\) vertices is occupied with
probability \(\mu^{-k}\). These are exactly the intersection
probabilities above. This is the one-dependent hard-core law on this
graph; its uniqueness from the one-site marginals is a general fact
\cite[Lemma~23]{HL}.

Another event has a simple interpretation without choosing a parity.
Call an interval of \(2j+1\) consecutive labels \emph{exposed} if
one upper and one lower arch leave the interval. For prescribed
disjoint intervals \(I_a\) of lengths \(2j_a+1\), the exact probability is
\begin{equation}\label{eq:blocks}
 \Prb(\text{all }I_a\text{ exposed})
    =\left(\prod_a M_{j_a+1}\right)\frac{M_{n-\sum_a j_a}}{M_n}.
\end{equation}
Indeed the interior of each interval is a path: an interior cycle would
already be a full connected component. Its two free incidences are
exposed to their respective half-planes. Adding a single vertex just
before the interval and closing these incidences yields an order-\((j_a+1)\)
meander. Deleting the added vertex inverts this operation. Contract
the interiors independently, retaining one vertex per interval, to
obtain~\eqref{eq:blocks}. Thus fixed disjoint intervals of this kind
are asymptotically independent, with marginal limits \(M_{j+1}/\mu^j\).

\section{The two orientation counts}
\label{sec:joint}
Write \(X_n^+,X_n^-\) for the positive and negative counts. For \(u,v\ge0\) let
\[
 Z_n(u,v)=\sum_{P\in\M_n}u^{X_n^+(P)}v^{X_n^-(P)},\qquad
 f_{u,v}(z)=z\left(\frac1{1-(u-1)z}+
                        \frac1{1-(v-1)z}-1\right).
\]
Here \(\M_n\) is the set of connected diagrams, and \(0^0=1\).
Let \(R=1/\mu\) and \(\M(z)=\sum_{n\ge1}M_nz^n\).

\begin{theorem}\label{thm:joint}
As a formal power series,
\[
             \sum_{n\ge1} Z_n(u,v)z^n=\M(f_{u,v}(z)).
\]
For \(u,v>0\), let \(\rho(u,v)\) be the first positive solution of
\(f_{u,v}(\rho)=R\), before any positive pole. Then
\(\lim Z_n(u,v)^{1/n}=1/\rho(u,v)\). Under the uniform measure,
\((X_n^+/n,X_n^-/n)\) satisfies a large-deviation principle at speed \(n\)
with good rate function
\[
 I(a,b)=\sup_{s,t\in\mathbb R}
       \bigl\{as+bt+\log\bigl(\mu\,\rho(e^s,e^t)\bigr)\bigr\}.
\]
\end{theorem}
\begin{proof}
Give each marked positive site weight \(u-1\), each marked negative
site weight \(v-1\), and sum over subsets of actual occurrences.
On contracting the marked runs, every even vertex of a diagram of
order \(m\) has three possibilities: it remains unchanged, carries a
nonempty positive run, or carries a nonempty negative run. A run of
length \(j\) increases the order by \(j\). The \(m\) vertices can be
expanded independently by Lemma~\ref{lem:contract}, which gives the
formal identity with both activities.

In the upper-arch splice of two diagrams the only positive occurrence
which can be lost is the first site of the right diagram. A negative
occurrence can be gained only at the seam. Indeed every short deleted
upper arch starts at an even label, the long new arch has length at
least three, and the new short upper arch starts at an odd label.
The lower matching remains unchanged. Removing the two crossing upper
arches inverts the splice at its specified position. Thus, for
\(u,v>0\),
\[
 Z_{r+s}(u,v)\ \geq\ \min(1,u^{-1})\min(1,v)\,Z_r(u,v)Z_s(u,v).
\]
Fekete's lemma, applied after this constant rescaling, gives root
convergence. The counts have an exponential upper bound by the
Catalan-square estimate and the bounded number of sites.

For nonnegative activities and real \(z\) before the first pole,
\[
 f'_{u,v}(z)=\frac1{(1-(u-1)z)^2}+\frac1{(1-(v-1)z)^2}-1
             \ \geq\ 2/(1+z)^2-1.
\]
On \([0,\sqrt2-1)\) the right-hand side is positive. If a positive pole
occurs within this interval the function tends to infinity there.
Otherwise \(f_{u,v}\ge f_{0,0}=z(1-z)/(1+z)\) and the latter reaches
\(3-2\sqrt2>R\) at its maximum, since \(\mu>11\) \cite{AP}. This proves
existence and uniqueness of the first root on the increasing branch,
with strictly positive derivative. The composition is analytic along
the positive interval before this root. Its coefficients are
nonnegative, so Pringsheim's theorem rules out a smaller radius.
Analyticity across the root would, by the local inverse, continue
\(\M\) across its positive singularity at \(R\). Hence the radius is
exactly \(\rho\). This radius argument includes zero activities.

The limiting logarithmic moment generating function for the two
counts is \(-\log(\mu\rho(e^s,e^t))\). It is differentiable on all of
\(\mathbb R^2\) by the implicit function theorem. The rescaled counts
lie in a compact triangle. The G\"artner--Ellis theorem
\cite[Theorem 2.3.6]{DZ} now gives the asserted full large-deviation
principle for open and closed sets.
\end{proof}

\begin{corollary}\label{cor:double}
Let \(A_n=Z_n(0,0)\) count diagrams avoiding both signs. Then
\[
 \lim_{n\to\infty} A_n^{1/n}
       = \frac{\mu-1+\sqrt{\mu^2-6\mu+1}}{2}.
\]
\end{corollary}
\begin{proof}
In a connected left factor of order \(r\ge2\), its last pair of
vertices cannot have both an upper and a lower short arch: that would
be an isolated cycle. If the lower arch is absent, use the upper
splice, so no negative site appears at the seam. Otherwise use the
lower splice; the absence of the last upper arch prevents a positive
site at the seam. The chosen colour is recovered from the two arches
crossing the cut. Thus \(A_{r+s}\ge A_r A_s\) for \(r\ge2,s\ge1\).
Here \(A_1=1,A_2=0,A_3=2\), so every \(A_n\) for \(n\ge3\) is
positive. For fixed \(k\ge3\) write \(n=jk+\ell\) with
\(3\le\ell\le k+2\); this proves root convergence by the usual
supermultiplicative argument. Its limit is the reciprocal of the
first solution of \(z(1-z)/(1+z)=1/\mu\).
\end{proof}

\section{Binomial shifts of the meander numbers}\label{sec:idle-moves}

Write \(M(z)=\sum_{n\ge1}M_nz^n\), where \(2n\) is the number of
crossings, and put
\[
     F_t(z)=M\!\left(\frac{z}{1+tz}\right),\qquad
     F_{t,n}=[z^n]F_t(z),\qquad C_n=F_{2,n}.
\]
Thus \(F_{0,n}=M_n\), and the one-sign avoidance identity gives
\(F_{1,n}=B_n\).

\begin{theorem}\label{thm:idle-shifts}
All coefficients of \(F_t\) are nonnegative if and only if \(t\le2\).
For each real \(t\le2\) they are eventually positive and
\[
       \lim_{n\to\infty}\frac{F_{t,n+1}}{F_{t,n}}=\mu-t .
\]
In particular the assertion at \(t=2\) holds despite \(C_2=0\).
Writing \(E M_n=M_{n+1}\), one has
\((E-2)^kM_n\ge0\) for every \(k\ge0\) and \(n\ge1\).
\end{theorem}

\begin{proof}
Use the four-letter sweep of a closed meander
\cite{Jensen,AP}. Turn the fixed line vertical and sweep upwards.
The state records the ordered ends on each side of the line and which
two ends belong to the same unfinished arc. A letter \(U\) creates a
fresh arc with one end on either side. A letter \(D\) joins the two
nearest ends on opposite sides. The letters \(R,L\) move the nearest
end from left to right and from right to left, respectively.
An unfinished word may never close a component. The first letter is
\(U\), giving a state \(e\) with one arc. The last letter is \(D\),
and its preceding state must be \(e\). The choices of letters determine
the two noncrossing matchings, since an end always closes against the
nearest available end on its side.

After an odd number of crossings there is at least one end on either
side. Indeed, on each side the number of ends equals the number of
crossings processed minus twice the number of completed matching edges
on that side. Divide the \(2n-2\) interior letters into consecutive
disjoint pairs, beginning immediately after the first \(U\). Write
\(T_{xy}\) for the number of legal two-letter words taking odd-cut
state \(x\) to state \(y\). Each row has sum at most sixteen, so every
entry of every power used below is a finite sum. Bijectivity of the
sweep gives \(M_{m+1}=(T^m)_{ee}\). At each odd cut the distinct pairs
\(LR\) and \(RL\) are legal and return the original state. Delete these
two choices at every state and call the resulting matrix \(Q=T-2I\).
In particular \(d_m=(Q^m)_{ee}\) counts meanders of order \(m+1\)
with neither pair in the designated positions.

The binomial theorem now gives
\[
   M_{m+1}=\sum_{j=0}^m\binom mj 2^{m-j}d_j,\qquad
   (E-2)^kM_n=(T^{n-1}Q^k)_{ee}.
\]
The first identity gives \(d_m=C_{m+1}\) by binomial inversion, and
the second proves the finite-difference assertion. More generally, for
\(s=2-t\) another binomial expansion gives
\begin{equation}\label{eq:idle-weight}
   F_{t,m+1}=\big((Q+sI)^m\big)_{ee}
            =\sum_{j=0}^m\binom mj s^{m-j}d_j .
\end{equation}
This follows also by substituting \(z/(1-sz)\) in \(F_2\).
It proves nonnegativity when \(s\ge0\). The coefficient
\(F_{t,2}=M_2-tM_1=2-t\) gives the converse.

Concatenating two return words at \(e\) gives
\(d_{m+h}\ge d_md_h\). The first four values are
\(d_0=1,\ d_1=0,\ d_2=4,\ d_3=10\), by binomial inversion from
\(M_1=1,\ M_2=2,\ M_3=8,\ M_4=42\). Every integer \(m\ge2\) is
a sum of twos and threes, so \(d_m>0\) in this range. Put
\(c=\sup_{k\ge2} d_k^{1/k}\le16\). For fixed \(k\ge2\) and all
large \(m\), write \(m=ak+b\) with \(2\le b\le k+1\).
Then \(d_m\ge d_k^a d_b\), so \(d_m^{1/m}\to c\) and \(c\ge2\).
For \(s>0\), the \(m\)-th root of the right side of
\eqref{eq:idle-weight} tends to \(s+c\). The upper bound follows
from \(d_j\le c^j\). For the lower bound use \(d_j\ge(c-\epsilon)^j\)
for all sufficiently large \(j\), with \(0<\epsilon<c\). The finitely
many omitted terms are \(O(m^K s^m)\) for fixed \(K\), which has
smaller exponential growth than \((s+c-\epsilon)^m\). Setting \(s=2\)
identifies \(c=\mu-2\). For any \(s>0\), the sequence
\(a_m=s^{-m}F_{t,m+1}\) is the binomial transform of the nonnegative
sequence \(s^{-j}d_j\) with initial term one, and its root tends to
\(1+c/s>1\). Lemma~\ref{lem:binomial-smoothing} gives
\(F_{t,m+2}/F_{t,m+1}\to s+c=\mu-t\).

It remains to treat \(s=0\). Two longer return words are available
at every odd cut even after the two original pairs have been deleted:
\[
               UR\mid DL,\qquad UU\mid RD\mid DL .
\]
Here the bars mark the disjoint pairs. To verify both words, note that
\(URD\) has exactly the effect of \(R\) whenever a left end is present.
The new arc created by \(U\) is separate from every old component.
The next letter moves its left end to the right. The letter \(D\)
joins that end to the original nearest left end, extending the old
arc by the fresh arc and placing its end on the right. The order and
pairing of the other ends are unchanged. Thus \(URDL\) returns the
original state, as in Albert and Paterson's zero-displacement example
\cite{AP}. In \(U(URD)DL\) the inner \(URD\) acts on the fresh
left end of the outer \(U\), and replacing it by \(R\) leaves \(URDL\).
Both constructions are legal even at \(e\). In particular
\(Q^r-I\) is entrywise nonnegative for \(r=2,3\).

Fix \(r\in\{2,3\}\) and use \(j=2,\ldots,r+1\) as representatives
of all residue classes modulo \(r\). Then
\[
 d_{rm+j}
    =\sum_{\ell=0}^m\binom m\ell
         \big(Q^j(Q^r-I)^\ell\big)_{ee}.
\]
The coefficients are nonnegative, their initial term is \(d_j>0\),
and the \(m\)-th root of the left side tends to \(c^r>1\).
Lemma~\ref{lem:binomial-smoothing} therefore gives
\(d_{n+r}/d_n\to c^r\) along every residue class and hence along
all integers. Finally,
\[
    \frac{d_{n+1}}{d_n}
       =\frac{d_{n+3}/d_n}{d_{n+3}/d_{n+1}}\longrightarrow c .
\]
Since \(d_n=C_{n+1}\), this proves the endpoint assertion.
\end{proof}

\begin{corollary}\label{cor:minimum-return}
For every \(k\ge0\), the minimum diagonal entry of \(Q^k\) equals
\((Q^k)_{ee}=C_{k+1}\).
\end{corollary}
\begin{proof}
Take a word returning \(e\) to itself. At a general odd cut mark the
nearest end on either side as an input. At each prefix the word uses only
these ends and their successors, since in the run from \(e\) it always
finds an available end on the required side. Its diagram in a rectangle
therefore has two inputs and two outputs. It has no internal circle.
Joining its inputs would make a circle when the starting arc of \(e\)
is restored, so both inputs pass through. Planarity fixes the order of
the two outputs. Gluing this rectangle to any original signature recreates
that signature. No prefix can close a component, since a closed component
persists. Every word counted at \(e\) is therefore counted on each
diagonal; the entry at \(e\) attains the bound.
\end{proof}

There is also a direct connection with ratio limits for countable
Markov kernels. Restrict the odd-cut transfer \(T\) of
Theorem~\ref{thm:idle-shifts} to the states communicating with \(e\);
every return from \(e\) stays in this class. Then \(T/16\) is
substochastic with diagonal at least \(2/16\). Kesten's lemma gives
the ratio limit for the ordinary meander numbers immediately
\cite[Lemma~4]{K95}; see also \cite[Lemma~1.2]{MSM}.
Deleting only the two specified loops preserves the communicating
class. The returns of lengths two and three in the proof show that
\(Q/16\) is uniformly aperiodic as well, giving the endpoint
ratio by the same lemma. The proof above makes the smoothing
argument explicit without using a general Markov-chain theorem.

\subsection*{Consecutive ratios at all nonnegative activities}
\begin{theorem}\label{thm:all-weighted-ratios}
For \(u,v\ge0\), let \(\rho(u,v)\) be the first positive solution of
\(f_{u,v}(\rho)=1/\mu\), before any positive pole. Then
\[
       \lim_{n\to\infty}\frac{Z_{n+1}(u,v)}{Z_n(u,v)}
                 =\frac{1}{\rho(u,v)} .
\]
In particular the ratio for avoidance of both signs tends to
\((\mu-1+\sqrt{\mu^2-6\mu+1})/2\).
\end{theorem}
\begin{proof}
First consider the root limit. Theorem~\ref{thm:joint} and
Corollary~\ref{cor:double} give it in the interior and at the origin.
On the axis \(u=0,\ v>0\), the upper splice used in the proof of
Theorem~\ref{thm:joint} preserves positive avoidance. It can add one
negative occurrence and loses none. Thus
\(Z_{r+s}(0,v)\ge\min(1,v)Z_r(0,v)Z_s(0,v)\).
Reflection gives the other axis. Fekete's lemma and the radius
calculation in that theorem give the root limit \(g=1/\rho(u,v)\)
throughout. Since \(Z_n(u,v)\ge A_n\), this limit exceeds one and
all sufficiently large terms are positive.

Use the four crossing moves of Theorem~\ref{thm:idle-shifts},
from \href{https://arxiv.org/pdf/cond-mat/0008178#page=7}{Jensen's
transfer} and the word description of Albert and Paterson \cite{AP}.
In a signature, list the exposed strand ends across the cut, with a
divider indicating the line. The letter \(U\) inserts the two ends of a
new strand at the divider. The letter \(D\) joins the two ends next to
the divider, provided they belong to distinct strands. The letter \(R\)
moves the last left end to the right and \(L\) moves the first right
end to the left. Pair the internal letters of a meander word starting
immediately after its initial \(U\), omitting its final \(D\).
Every tested site is a triple \(pab\), where \(ab\) is one of these
pairs and \(p\) is the preceding letter. One sign occurs precisely when
\[
         p\in\{U,L\},\quad a=R,\quad b\in\{D,L\},
\]
and the other when \(p\in\{U,R\},\ a=L,\ b\in\{D,R\}\).
Indeed the two short arches require an opening and a closing in
consecutive positions, first on one side and then on the other.
Identify the left side of the list with the upper half-plane of
the arch diagram. Then \(R\) closes above and opens below, so the
first test has weight \(u\), the second \(v\). Reflection exchanges
the tests.

Record the previous letter together with the signature. Let \(W\)
count allowed pairs on this state space, multiplying each pair by
the activity of its sign, or by one when neither sign occurs. At a
pair boundary the number of ends on each side is odd. The initial
state is a single strand across the divider with previous letter
\(U\). It cannot occur after a pair, since any subsequent \(U\)
inserts a new strand in a nonempty signature. Let \(\alpha\) be the
weighted row after the first pair from this initial state, and let
\(\beta\) indicate the single-strand signature across the divider
with any attainable previous letter. Jensen's final closure then gives
\[
                  Z_{N+2}(u,v)=\alpha W^N\beta \qquad (N\ge0).
\]
Here \(W\) acts on the states reachable after at least one pair,
with all original strand pairings retained. All products are
entrywise finite by the four-letter rule.

For every such state we exhibit two returns that preserve its previous
letter and create no marked site. In the table, a bar divides the
exposed ends. The letters \(a,b\) represent ends of pre-existing
strands. Equal digits are the two ends of one local strand. The names
of local strands can be reassigned after each pair. The last column
records the exposed ends after each successive pair of letters.
\begin{center}
\small
\begin{tabular}{@{}llll@{}}
\toprule
Previous & Initial ends & Pairs & Ends after the pairs\\
\midrule
\(\mathrm R\)&\(\mathrm{a|b}\)&
\(\mathrm{UL,DR}\)&\(\mathrm{a11|b,\ a|b}\)\\
& &\(\mathrm{UL,LR,DR}\)&\(\mathrm{a11|b,\ a11|b,\ a|b}\)\\[2pt]
\(\mathrm D\)&\(\mathrm{a|b}\)&
\(\mathrm{LU,UR,DD}\)&\(\mathrm{ab1|1,\ ab1|221,\ a|b}\)\\
& &\(\mathrm{UU,LL,DR,RD}\)&\(\mathrm{a12|21b,\ a1221|b,\ ab1|1,\ a|b}\)\\[2pt]
\(\mathrm U\)&\(\mathrm{ab1|1}\)&
\(\mathrm{UR,DD,LU}\)&\(\mathrm{ab1|221,\ a|b,\ ab1|1}\)\\
& &\(\mathrm{UL,DR,RD,LU}\)&\(\mathrm{ab122|1,\ ab1|1,\ a|b,\ ab1|1}\)\\
\bottomrule
\end{tabular}
\end{center}
Reflection gives the two returns for previous letter \(L\) from the
\(R\) rows. When the previous letter is \(U\), the central ends are
paired. There is at least one older strand, and each side has an even
number of other ends. Hence one side has at least two. Display the
two nearest on that side as \(a,b\), reflecting the \(U\) rows if
necessary. Every other exposed end stays outside the calculation.
For previous letter \(R,L\) or \(D\), a single nearest end on each
side suffices. No step joins \(a\) to \(b\), so \(a\) and \(b\) may
also be paired in the full signature. Each
listed word ends in its prescribed previous letter. Checking its pairs
in order against the two triples above gives weight one, including
the first pair at the incoming boundary.

The return lengths in pairs are \(2,3\) or \(3,4\). Repeating a
return of length 2 gives length 4, so \(W^d\ge I\) entrywise for
\(d=3,4\), for every \(u,v\ge0\).
Fix either \(d\), and for each residue class take a sufficiently large
representative \(r\) that \(\alpha W^r\beta>0\). Then
\[
 \alpha W^{r+dj}\beta
   =\sum_{i=0}^j\binom ji\alpha W^r(W^d-I)^i\beta .
\]
The binomial weights on the right are nonnegative, the zeroth is
positive, and the \(j\)-th root of the left side tends to \(g^d>1\).
Lemma~\ref{lem:binomial-smoothing} gives
\(Z_{n+d}/Z_n\to g^d\) on each residue class and hence on their
finite union. Finally
\[
      \frac{Z_{n+1}}{Z_n}
        =\frac{Z_{n+4}/Z_n}{Z_{n+4}/Z_{n+1}}\ \longrightarrow\ g .
\]
\end{proof}

\begin{corollary}\label{cor:uniform-weighted-ratios}
The convergence in Theorem~\ref{thm:all-weighted-ratios} is uniform on
every compact subset of \([0,\infty)^2\), including sets meeting either
axis.
\end{corollary}
\begin{proof}
Write \(\theta=(u,v)\) and \(g(\theta)=1/\rho(u,v)\).
The first root of \(f_{u,v}(z)=1/\mu\) is simple, lies on the
increasing branch, and precedes the first positive pole, as shown in
Theorem~\ref{thm:joint}. The implicit function theorem therefore gives
continuity of \(g\) up to the boundary of the quadrant. For every
\(1<c<g(0,0)\) and every sufficiently large \(n\),
\(Z_n(\theta)\ge A_n\ge c^n\) throughout the quadrant.
Moreover, if \(n_j\to\infty\) and \(\theta_j\to\theta\), then
\[
       Z_{n_j}(\theta_j)^{1/n_j}\longrightarrow g(\theta).
\]
Indeed the coefficients of \(Z_n\) are nonnegative. Bound
\(\theta_j\) coordinatewise between
\(\theta^-=\max(0,\theta-\delta\mathbf 1)\) and
\(\theta^+=\theta+\delta\mathbf 1\), apply the root limit at those two
fixed points, then let \(\delta\downarrow0\), where \(\mathbf 1=(1,1)\).

Fix \(d\in\{3,4\}\). In each residue class modulo \(d\) choose
\(r\ge1\) with \(A_{r+2}>0\). In the transfer notation of the theorem put
\[
 a_j(\theta)=Z_{r+dj+2}(\theta),\qquad
 b_i(\theta)=\alpha W^r(W^d-I)^i\beta,\qquad
 a_j=\sum_{i=0}^j\binom ji b_i .
\]
The unmarked return at every state gives a constant term at least one
in every diagonal entry of \(W^d\). Thus each \(b_i\) is a polynomial
with nonnegative coefficients, and
\(b_0(\theta)\ge A_{r+2}>0\). Along any \(\theta_j\to\theta\) we have
\(a_j(\theta_j)^{1/j}\to\lambda=g(\theta)^d>1\).

Give \(J_j=i\) probability \(\binom ji b_i(\theta_j)/a_j(\theta_j)\).
At the fixed upper point \(\theta^+\), the generating-function
identity in Lemma~\ref{lem:binomial-smoothing} gives
\(\limsup_i b_i(\theta^+)^{1/i}=g(\theta^+)^d-1\).
Hence for any \(\epsilon>0\), all sufficiently large \(j\) and all \(i\),
\[
 b_i(\theta_j)\le b_i(\theta^+)
            \le C_{\delta,\epsilon}\bigl(g(\theta^+)^d-1+\epsilon\bigr)^i .
\]
Write \(H(t)=-t\log t-(1-t)\log(1-t)\) with \(0\log0=0\).
For any nonempty closed \(F\subset[0,1]\), the entropy bound on
\(\binom ji\) and the root limit for \(a_j(\theta_j)\) imply
\[
 \limsup_{j\to\infty}j^{-1}\log\Pr(J_j/j\in F)
 \le \max_{t\in F}\{H(t)+t\log(\lambda-1)-\log\lambda\}.
\]
Here one first uses \(g(\theta^+)^d-1+\epsilon\) in place of
\(\lambda-1\), then sends \(\delta,\epsilon\) to zero. The continuous
expression on the right has unique maximum zero at \(1-1/\lambda\) on
\([0,1]\). Thus \(J_j/j\) concentrates there and its expectation
converges to that point. The binomial ratio identity gives
\[
       \frac{a_{j-1}(\theta_j)}{a_j(\theta_j)}
                =1-\mathbb E[J_j/j]\longrightarrow g(\theta)^{-d}.
\]
Every sequence of parameters in a compact set admits a convergent
subsequence, and the residue classes are finite. Since \(g\) is
continuous and finite there, this proves uniform convergence of the
forward \(d\)-step ratios to \(g^d\). Taking the quotient of the
four-step ratio at \(n\) and the three-step ratio at \(n+1\) gives
the assertion, since the latter tends uniformly to \(g^3\ge g(0,0)^3>1\).
\end{proof}

\begin{corollary}
Fix \(d=3\) or \(4\). For each attainable state after the first pair,
choose one of the \(d\)-pair unmarked returns constructed in the proof
(for length four, a two-pair return may be repeated). At \(k\)
prescribed disjoint intervals of \(d\) pairs, all after the first
pair, the probability of seeing the chosen word at every interval
under the \((u,v)\)-weighted law is exactly
\[
                    \frac{Z_{n-dk}(u,v)}{Z_n(u,v)}
\]
whenever the placement is possible and the denominator is positive.
For fixed \(k\) this tends to \(\rho(u,v)^{dk}\).
\end{corollary}
\begin{proof}
Delete the indicated words. Each restores its full incoming state,
including the previous letter, and contributes weight one.
Conversely the prescribed positions and the state there determine
each inserted word uniquely. This bijection preserves the weight of
the rest of the diagram. The ratio limit gives the assertion. Applying
the same formula to each subset shows joint convergence of the indicators
on disjoint intervals to independent Bernoulli variables of parameter
\(\rho(u,v)^d\).
\end{proof}

\section{Conditioning on avoidance of one orientation}

Let \(B_n\) count connected meanders with no positive site, and put
\(b=\mu-1\), so \(B_n^{1/n}\to b\). Give this avoidance class the
uniform probability measure. Write \(Y\) for the number of negative
sites under this measure. In this section the site set is
\(V_n=\{0,\ldots,n-2\}\).

\begin{theorem}\label{thm:conditional-ratio}
The one-orientation avoidance counts satisfy
\[
             \lim_{n\to\infty}\frac{B_n}{B_{n-1}}=b .
\]
In the interior, the negative-site indicators converge to a stationary
two-dependent process. For a finite set \(T\subset\mathbb Z\), with
\(k=|T|\) and \(d=|((T-1)\cup(T+1))\setminus T|\), its intersection
probabilities are
\[
   \mathbb P(\text{all sites of }T\text{ are negative})
             =b^{-k}(1+b^{-1})^d .
\]
In particular its one-site density equals \((b+1)^2/b^3\).
\end{theorem}

\begin{proof}
At finite order fix \(T\subseteq V_n\) and take the exterior neighbours
only in \(V_n\), so that
\(d=|(((T-1)\cup(T+1))\cap V_n)\setminus T|\).
A negative site excludes a positive site at its own location and at
each neighbour. There remain \(L=n-1-k-d\) sites at which a positive
occurrence must be excluded. Fixed-site contraction and
inclusion--exclusion give the count
\[
 \sum_{j=0}^{L}(-1)^j\binom Lj M_{n-k-j}
       =\sum_{a=0}^{d}\binom da B_{n-k-a}.                 \tag{*}
\]
For the equality, the left-hand sum also counts order-\((n-k)\)
meanders avoiding positives at any prescribed \(L\) sites.
Their exact positive set is a subset of the other \(d\) sites.
The exact-set formula gives the right-hand sum.

The ratio follows from Theorem~\ref{thm:idle-shifts} at \(t=1\).
Divide \((*)\) by \(B_n\) and use the consecutive-ratio limit for
each of its finitely many terms. Sets \(T\) translated into the
interior retain both exterior neighbours whenever these are distinct.
The resulting intersection probabilities determine every finite law by
inclusion--exclusion and are translation invariant. If two sets have
mutual distance at least three, both \(k\) and \(d\) add for each choice
of subsets; their joint law factors. Singletons at distance two have
three distinct exterior neighbours in total, giving strictly less than
the product of their one-site probabilities.
\end{proof}

The same law arises from conditioning in a finite window of the
unconditioned limiting coloured process. Let \(p=1/\mu\), let \(T\) be
finite, and let \(W\) contain \(T\) and its neighbours. Prescribing
negatives on \(T\) excludes positives on \(k+d\) sites of \(W\).
Inclusion--exclusion on the remaining sites gives the intersection
probability \(p^k(1-p)^{|W|-k-d}\). Dividing by
\(\mathbb P(\text{no positive in }W)=(1-p)^{|W|}\) gives exactly
\(b^{-k}(1+b^{-1})^d\), independently of how far \(W\) extends beyond
the neighbours. Thus conditioning the finite meander on avoidance
throughout its full site set has the same local limit as this
finite-window calculation.

\subsection*{A finite entropy comparison in the full class}
For \(n\ge2\) put \(N=n-1\), \(p=1/\mu\), and \(q_n=M_{n-1}/M_n\).
Writing \(D\) for binary relative entropy, one has
\[
 N D(q_n\Vert p)\le
      \log\!\left(\frac{b\mu^{\,N}}{M_n}\right),\qquad
 |q_n-p|\le
      \sqrt{\frac{1}{2N}\log\!\left(\frac{b\mu^{\,N}}{M_n}\right)}.
\]
Indeed \(B_j\le b^j\) by supermultiplicativity. A particular binary
indicator word of weight \(k\) has probability \(B_{n-k}/M_n\),
whereas its probability under the product Bernoulli-\(p\) law is
\(b^{N-k}/\mu^N\). Thus the relative entropy of the full law with
respect to the product law is at most the logarithm displayed above.
Subadditivity of Shannon entropy bounds it below by the sum of the
\(N\) one-site divergences. Finally \(D(q\Vert p)\ge2(q-p)^2\),
as follows by twice differentiating in \(q\). The right-hand side of
the first inequality is \(o(n)\) by root growth. This also proves the
unconditioned ratio theorem directly. A numerical error bound from
this expression requires information about \(\mu\).

\subsection*{The longest run of one orientation}
Set \(N=n-1\), and let \(L_n\) be the longest run of positive sites.
Conditional on \(X_n^+=x\), the indicator word is uniform among the
\(\binom Nx\) words with \(x\) ones. The lengths of the runs of ones,
including zero-length gaps, form a uniform weak composition of \(x\)
into \(g=N-x+1\) parts. If \(Y_r\) counts parts of length at least
\(r\ge1\), subtract \(r\) from each of \(j\) ordered selected parts:
\[
 {\mathbb E}[(Y_r)_j\mid X_n^+=x]=(g)_j\frac{(x)_{jr}}{(N)_{jr}}
 \quad(j\le g,\ jr\le x).
\]
Here the parentheses denote falling factorials; outside this range
the moment is zero. For \(x/N\) in a compact subinterval of \((0,1)\)
and \(r=O(\log N)\), the right-hand side is
\([g(x/N)^r]^j(1+O_j((r+1)^2/N))\) for fixed \(j\).
At \(r\) just above \((1+\epsilon)\log N/\log(N/x)\), the first
moment tends to zero. Just below
\((1-\epsilon)\log N/\log(N/x)\), the first moment diverges and the
ratio of the second factorial moment to its square tends to one.
Rounding to integers, Markov's inequality and the second moment
bound place \(L_n\) between these thresholds with probability tending
to one. Since \(X_n^+/N\to1/\mu\), we obtain
\[
                         L_n/\log n\ \longrightarrow\ 1/\log\mu
                         \qquad\hbox{in probability}.
\]

\section{Finite tests and related work}
An exhaustive ordered-matching calculation checks the two-dimensional
histogram for every \(1\le n\le8\). The largest catalogue contains
\(1{,}430^2=2{,}044{,}900\) pairs before the connectivity test.
At order eight there are \(n-1=7\) possible sites; the 128 prescribed
subsets for the conditional law were checked individually.
A separate inverse-splice implementation tests all admissible seam
positions in the small catalogues. In the supplied finite tables the
mixed binomial moments are obtained by distributing \(i\) positive
and \(j\) negative marks in nonempty runs among the \(n-i-j\)
remaining vertices. Two integer binomial inversions give the full
joint histogram. The tables through order 38 use the existing
values of \(M_n\) \cite{OEIS}.

A second transfer records the preceding letter and both mark counts.
Its coefficients agree with the joint table through order eleven.
The supplement tests the unmarked returns of Theorem~\ref{thm:all-weighted-ratios}
on every noncrossing matching with up to sixteen boundary ports and each
appropriate odd split. To fix the sign convention, the code also compares
the marks of individual ordered pairs with their sweep words through
order seven. A separate check tests the 741 shifted finite differences
available from the ordinary counts through order 38.

Lando and Zvonkin \cite[Sec.~6]{LZ} studied the successive ratios
numerically by rational interpolation. Albert and Paterson \cite{AP}
use local moves and cluster enumeration to bound the growth constant.
Caroppo, Da Lozzo, Di Battista and Vihrovs (\href{https://arxiv.org/abs/2610.11944}{arXiv:2610.11944}, 2026) extend that calculation with longer
representatives of forbidden words.
Belousov's prime decomposition for open meanders includes inserting
odd inverse snakes at several intersections to bound the growth of
almost irreducible objects \cite[Thm.~8]{Bel}. Delecroix et al.\
\cite{DGZZ} study an asymptotic regime with a fixed number of minimal
arches. Féray and Thévenin \cite{FT} obtain local convergence for
uniform meandric systems sampled without conditioning on connectivity.
Slote's operadic construction (\href{https://joeslote.com/documents/operadic_techniques_meanders_WIP.pdf}{\emph{Operadic Techniques for Meander Enumeration}}, 2019, Definition~2) treats
insertion in a chosen input of a closed meander. The marked-run
series is related to classical cluster marking of overlapping
patterns \cite{GJ}. The use of local patterns in ratio-limit arguments has precedents in
lattice models, notably Kesten's work on self-avoiding walks \cite{K}
and Madras's pattern theorem for clusters \cite{Mad}. Here the
fixed-site contraction gives an exact law under connectivity conditioning.
The argument determines neither a numerical value of \(\mu\) nor the
power correction to exponential growth.

\section*{Author statements}
\paragraph{Funding.}
The author funded this work personally.
\paragraph{Competing interests.}
The author declares no competing interests.
\paragraph{Author contribution.}
Iaroslav Argunov is the sole author. He conducted the research and checked
the proofs and references. Computational and writing assistance is described below.
\paragraph{Data and software.}
The joint count tables, the Python verification programs and the reference
outputs for the finite checks are available in the
\href{https://github.com/Kimargin8/meander-ratio-and-smith/tree/72c5398f8c1e9cf5966fcbbf558f34737497bc01/ratio}{versioned accompanying repository}.
The ordinary meander counts used as input come from the On-Line
Encyclopedia of Integer Sequences, \href{https://oeis.org/A005315}{A005315},
orders 1 through 38.
\section*{Acknowledgements}
The author used OpenAI Codex to help develop proof arguments, write and
check the enumeration and verification programs, search the literature,
and draft and edit the English manuscript. The author is responsible for
the arguments, references and text.
\end{document}